\documentclass[11pt]{amsart}

\usepackage{amssymb, upgreek, verbatim, graphicx, lmodern} 
\usepackage[left=1in,top=1in,right=1in,bottom=1in,head=.2in]{geometry}

\usepackage{enumitem}
\usepackage[textsize=scriptsize]{todonotes}
\usepackage{color}
\usepackage{mathtools} 
\usepackage{marginnote}

\usepackage{fancyhdr}
\usepackage{mathrsfs}
\usepackage[cmtip,all]{xy}
\usepackage{hyperref}

\swapnumbers
\newtheorem{theorem}{Theorem}[section] 

\newtheorem{corollary}[theorem]{Corollary}
\newtheorem*{theorem*}{Theorem}

\newtheorem*{fcthm*}{Finite Cork Theorem}
\newtheorem*{ccthm*}{Cork Consolidation Theorem}
\newtheorem*{csthm*}{Cork Separation Theorem}
\newtheorem*{icthm*}{Infinite Cork Theorem}
\newtheorem*{ucthm*}{Universal Cork Theorem}
\newtheorem*{aclemma*}{\ac\,Lemma}
\newtheorem*{smoothinglemma*}{Smoothing Lemma}
\newtheorem*{mclemma*}{Multicork Lemma}
\newtheorem*{multicorktheorem*}{Multicork Theorem}
\newtheorem*{lemma*}{Lemma}
\newtheorem*{corollary*}{Corollary}
\newcommand{\thistheoremname}{}
\newtheorem{genericthm}[theorem]{\thistheoremname}
\newenvironment{namedtheorem}[1]
  {\renewcommand{\thistheoremname}{#1}%
   \begin{genericthm}}
  {\end{genericthm}}

\theoremstyle{definition}
\newtheorem{definition}[theorem]{Definition}
\newtheorem{remark}[theorem]{Remark}
\newtheorem{remarks}[theorem]{Remarks}

\newtheorem*{remark*}{Remark}
\newtheorem*{definition*}{Definition}
\newtheorem*{remarks*}{Remarks}
\newtheorem*{addenda*}{Addenda}

\newcommand{\pf}{\vskip-5pt \vskip-5pt \proof}

\newcommand{\fig}[3]{\begin{figure}[h!] \includegraphics[height=#1pt]{#2}#3\end{figure}}
\newcommand{\defref}[1]{Definition~\ref{#1}}
\newcommand{\figref}[1]{Figure~\ref{#1}}

\newcommand{\thmref}[1]{Theorem~\ref{#1}}
\newcommand{\lemref}[1]{Lemma~\ref{#1}}

\newcommand{\corref}[1]{Corollary~\ref{#1}}
\newcommand{\remref}[1]{Remark~\ref{#1}}

\newcommand{\head}[1]{\bigskip\noindent{\bf #1}}
\newcommand{\headtop}[1]{\noindent{\bf #1}}

\newcommand{\bit}[1]{\textit{#1}} 

\newcommand{\bz}{\mathbb Z}

\newcommand{\calc}{\mathcal C}

\newcommand{\medcup}{\mbox{\larger$\cup$}}
\newcommand{\medsqcup}{\mbox{\larger$\sqcup$}}
\newcommand{\st}{\ \vert\ }
\newcommand{\id}{\textup{id}}

\newcommand{\lfrom}{\ \longleftarrow \ }
\newcommand{\lto}{\ \longrightarrow \ }
\newcommand{\hto}{\hookrightarrow}
\newcommand{\sto}{\!\!\xymatrix@C=1em{{}\ar@{~>}[r]&{}}\!\!}
\newcommand{\lhto}{\ \lhook\joinrel\relbar\joinrel\rightarrow \ }

\newcommand{\del}{\partial}
\newcommand{\dt}{{\hskip-.2pt\raisebox{.2ex}{\text{\scalebox{.6}{\textbullet}}}\,}}

\newcommand{\sss}{S^2\!\times\!S^2}
\newcommand{\sts}{S^2\widetilde\times S^2}

\newcommand{\interior}{\textup{int}}
\newcommand{\cl}{\textup{cl}}
\newcommand{\ac}{\textup{AC}}

\newcommand{\overbar}[1]{\mkern 1mu\overline{\mkern-2.5mu#1\mkern-1mu}\mkern 1mu}

\newcommand{\cs}{\mathop\#}
\newcommand{\bcs}{\mathop\natural}
\newcommand{\bsum}{C_0\bcs \cdots \bcs C_{n-1}}
\newcommand{\csum}{\del C_0 \cs \cdots \cs \del C_{n-1}}
\newcommand{\ddsum}{B^4\medsqcup C_0\medsqcup\cdots\medsqcup C_{n-1}}
\newcommand{\esum}{C_0\bcs B_1\bcs \cdots \bcs B_{n-1}}
\newcommand{\dualp}{{p^*}}

\DeclareMathOperator{\diff}{Diff}

\newcommand{\ia}{{\small\bf a)}}
\newcommand{\ib}{{\small\bf b)}}

\newcommand{\ione}{{\small I}}
\newcommand{\itwo}{{\small II}}

\def\X{{X^*}}
\def\Q{\sqcup\kern 0.1em C_i}
\def\B{\natural\kern 0.1em C_i}
\def\q{\medcup\kern 0.1em A_i}
\def\b{\natural\kern 0.1em A_i}
\def\Ki{\medcup K_i}
\def\Li{\medcup L_i}

\newcommand{\one}{X^{\scalebox{0.6}{\bf 1}}}
\newcommand{\dualone}{X^{*\scalebox{0.6}{\bf 1}}} 
\newcommand{\ks}{X^{\scalebox{0.6}{\bf k}}}
\newcommand{\dualks}{X^{*\scalebox{0.6}{\bf k}}}
\newcommand{\dri}{{r^*_i}}
\newcommand{\drj}{{r^*_j}}
\newcommand{\odri}{\overbar r^*_i}

\renewcommand{\H}{H^*}

\def\cl{\textup{c$\ell$}}

\newcommand{\foot}[1]{\setcounter{footnote}{1}\footnote{\ #1}}

\newcommand{\items}{\begin{itemize}[leftmargin=25pt,rightmargin=5pt]
  \setlength\itemsep{2pt}\vskip-3pt\vskip-3pt}
\newcommand{\stopitems}{\end{itemize}}

\begin{document}

\title{Higher Order Corks}
\author{Paul Melvin}
\address{Department of Mathematics\newline\indent  Bryn Mawr College\newline\indent  Bryn Mawr, PA 19010}

\author{Hannah Schwartz}
\address{Department of Mathematics\newline\indent  Princeton University \newline\indent  Princeton, NJ 08544}

\vspace*{-10pt}

\begin{abstract}
It is shown that any finite list of smooth closed simply-connected $4$-manifolds homeomorphic to a given one $X$ can be obtained by removing a single compact contractible submanifold (or \!{\it cork}) from $X$, and then regluing it by powers of a boundary diffeomorphism. We then use this result to `separate' finite families of corks embedded in a fixed $4$-manifold. \end{abstract}

\maketitle

\vskip-30pt
\vskip-10pt
\parskip 3pt

\setcounter{section}{-1}

\section{Introduction}  


A \bit{cork} is a compact contractible 4-manifold $C$ equipped with a boundary diffeomorphism $h\colon\del C\to\del C$.\foot{We work throughout in the smooth oriented category, so implicitly assume all manifolds are smooth and oriented, and all diffeomorphisms are orientation preserving.}   The cork $(C,h)$ is \bit{trivial} if $h$ extends to a diffeomorphism of $C$\,; for example $(B^4,h)$ is trivial for any $h$ \cite{cerf}.  The cork is \bit{finite} of order $n$ if $h$ is periodic of order $n$.  Corks of order $2$ will be called \bit{involutory}.      

If $C$ is embedded in the interior of a 4-manifold $X$, then the associated \bit{cork twist}
$$
X_{C,h} \ := \ (X-\interior(C))\,\cup_h \,C
$$
is homeomorphic \cite{freedman:simply-connected} but not generally diffeomorphic \cite{akbulut:contractible} to $X$.  If $X_{C,h}$ and $X$ {\it are} diffeomorphic, written $X_{C,h}\cong X$, then $C\subset X$ is said to be a \bit{trivial embedding} of the cork.  

It can be shown that a cork is nontrivial if and only if it embeds nontrivially in some $4$-manifold \cite{akbulut-ruberman:absolute}\cite{akmr:equivariant}.
The first nontrivial (involutory) cork was found by Akbulut \cite{akbulut:contractible}, showing that cork twists can alter smooth structures on 4-manifolds.  It is now known by Curtis-Freedman-Hsiang-Stong \cite{curtis-freedman-hsiang-stong} and Matveyev \cite{matveyev} that {\sl any} pair of compact simply-connected 4-manifolds that are h-cobordant rel boundary 
(so homeomorphic by Freedman \cite{freedman:simply-connected}) are related by a single involutory cork twist, and that for closed manifolds, the cork may be chosen with simply-connected complement.
We call this the ``Involutory Cork Theorem" (see \S\ref{ICT}).  

Our first main result extends this theorem to arbitrary finite lists of closed simply-connected $4\text{-manifolds}$ (or more generally compact simply-connected $4$-manifolds bounded by homology spheres, see \thmref{rfct}).  This demonstrates the ubiquity of nontrivial higher order corks, meaning corks $(C,h)$ of order $n\ge3$ whose twists $X_{C,h^1},\dots,X_{C,h^n}$  for suitable embeddings $C\subset X$ are smoothly distinct.  Such corks were first shown to exist in \cite{akmr:equivariant}; see Tange \cite{tange:finitecorks} for a weaker existence result. 

\begin{fcthm*}\label{fcthm}
For any finite list $X_i$ \,$(i\in\bz_n)$  of closed simply-connected  $4$-manifolds all homeomorphic to a given one $X = X_0$, there is a cork $(C,h)$ of order $n$ embedded in $X$ with simply-connected complement whose twists $X_{C,h^i}$ are diffeomorphic to $X_i$ for each $i\in\bz_n$.
\end{fcthm*}

Our main application of the Relative Involutory Cork Theorem \ref{ICT}, along with the Consolidation Theorem \ref{consolidation} (our main technical theorem from Section \ref{fcs}), is a recipe for finding \emph{disjoint} involutory corks in a closed, simply-connected $4$-manifold.

\begin{namedtheorem}{Separation Theorem}\label{separation} 
For any $n \in \mathbb{N}$ and any family of corks $(C_1,\sigma_1), \dots, (C_n,\sigma_n)$ embedded in a {\sl closed} simply-connected $4$-manifold $X$, there is a corresponding family of simple involutory corks $(A_i, \tau_i)$, embedded disjointly and each with simply-connected complement in $X$, whose twists $X_{A_i, \tau_i}$ are diffeomorphic to $ X_{C_i,\sigma_i}$ for each $i$. 
\end{namedtheorem}  

Although it is tempting to search for a single infinite cork in $X$ with cork twists diffeomorphic to an infinite list $X_i$ $(i\in\bz)$ of smooth $4$-manifolds homeomorphic to $X$ (such as an enumeration of {\sl all} the exotic smooth structures on $X$), such a cork need not exist, as noted by Tange \cite{tange}. For example, it follows from the adjunction inequality that knot surgeries on the Kummer surface using any list of knots with unbounded Alexander polynomial degrees cannot be the cork twists associated to a fixed embedding of a single infinite cork (cf.\ Yasui \cite{yasui}).  Thus the Finite Cork Theorem has no direct infinite analogue.

\vskip-4pt
\vskip-4pt
\head{Acknowledgments.}  This work was initiated while both authors were visitors at IAS in Princeton in the Fall of 2016.  We thank the Institute for its hospitality.  We are also indebted to Dave Auckly for enlightening conversations leading to the relative versions of our cork theorems, to Danny Ruberman and Bob Gompf for asking timely questions that inspired our separation theorem, and to the referee for a remarkably insightful and detailed report.

\section{Preliminaries}\label{prelims}

\headtop{\ac\ manifolds}

Let $X$ be a 4-manifold that can be built from the $4$-ball by adding $1$ and $2$-handles.  Any such handle structure is specified by a link $K\cup L$ in the $3$-sphere, where $K$ is a dotted unlink representing the $1$-handles\foot{The $1$-handles can be viewed as trivial $2$-handles carved out from the $4$-ball, see for example \cite[Chapter I.2]{kirby:4-manifolds}.  These $2$-handles are obtained by pushing into $B^4$ a family of disjoint spanning disks in $S^3$ for the dotted circles.  Such families of disks are {\sl not} unique (up to isotopy rel boundary) in $S^3$, but they are unique in $B^4$.} and $L$ is the framed link of attaching circles for the $2$-handles.  We write $X = [K,L]$ to specify this handle structure, which induces a presentation of the fundamental group
$$
\pi_1(X) \ = \ (\,x_1,\, \dots , \, x_k \st r_1,\,\dots, \, r_\ell\,)
$$  
where $x_1,\dots,x_k$ are the meridians of the components of $K$, and $r_1,\dots, r_\ell$ are the words in the free group $F(x_1,\dots,x_k)$ (determined up to conjugacy) traced out by the components of $L$.

\begin{definition}\label{acmanifold}
A handle structure $X=[K,L]$ as above is called an \bit{\ac\,structure} if some subset of the $2$-handles homotopically cancel some subset of the $1$-handles, while the remaining $2$-handles are attached homotopically trivially. That is, for some $n \leq \min(k,\ell)$ and a suitable ordering of the handles, $r_i = x_i$ for $i \le n$ and $r_i=1$ for $i>n$.

We call a $4$-manifold \bit{\ac}, or an \bit{\ac\,manifold}, if it admits an \ac\,structure.  Any such manifold is homotopy equivalent to a wedge of circles and $2$-spheres (namely $(\vee^{k-n}S^1)\vee(\vee^{\ell-n}S^2)$).  In the contractible case ($n=k=\ell$) we may sometimes refer to the manifold as an \bit{\ac\,cork}, even if it is not equipped with a boundary diffeomorphism.   Thus \ac\,corks are just corks (i.e.\ compact contractible $4\text{-manifolds}$) that can be built from the $4$-ball by adding an equal number of homotopically cancelling $1$ and $2$-handles.  It is unknown whether all corks are \ac; see the first remark below.  

For our purposes it is convenient to give a name to the two extreme {\it \ac\,structure types} when $n$ is equal to $\ell$ or to $k$, which overlap in the contractible ones:

{\bf Type \ione} \ ($n=\ell$) The $2$-handles homotopically cancel some subset of the $1$-handles, so $X\simeq\vee S^1$'s

{\bf Type \itwo} ($n=k$) Some subset of the $2$-handles homotopically cancel the $1$-handles, so $X\simeq\vee S^2$'s
\end{definition}

\begin{remarks} \label{andrewscurtis}
\ia\ The acronym \ac\ refers to the Andrews-Curtis moves on group presentations \cite{ac} that correspond to 1 and 2-handle slides: 
\items
\item[\dt] replace a generator by its inverse or by its product with another generator, and
\item[\dt] replace a relator by its inverse or by its product with a conjugate of another relator.
\stopitems
It is unknown whether all balanced presentations of the trivial group can be trivialized through Andrews-Curtis moves.  This is the {\it Andrews-Curtis Conjecture}.   Even the weaker {\it stable} version of this conjecture (which allows the addition of a new generator $x$ and relation $r=x$, corresponding to the introduction of a cancelling 1/2-handle pair) is unknown; it would imply that every cork built without $3$-handles (or without $1$-handles \cite{trace}) is $\ac$.

\ib\  The double $C\cup_\id -C$ of any \ac\,cork $C$ is diffeomorphic to $S^4$, as it is the boundary of $[0,1] \times C$, which is diffeomorphic to $B^5$ since homotopy implies isotopy for curves in 4-manifolds.  But there is no reason to suspect, unless the smooth $4$-dimensional Poincar\'e Conjecture holds, that all ``twisted doubles" $C\cup_{h}-C$ of $C$ for $h\in\diff(\del C)$ should be diffeomorphic to $S^4$.  Indeed every homotopy $4$-sphere can be constructed as a twisted double of some \ac\,cork $(C,h)$; see \corref{twisteddouble} below.
\end{remarks} 

\begin{definition}\label{simple}
A cork $(C,h)$ is \bit{simple} if it is \ac\ and $C\cup_{h^i}-C\cong S^4$ for all $i$.
\end{definition} 


\headtop{\ac\ cobordisms}

Now let $X$ be a cobordism between closed connected $3$-manifolds $M$ and $N$ that can be built from $M\times  [0,\varepsilon]$ by adding 1 and 2-handles along $M\times\{\varepsilon\}$.   Any such relative handle structure is specified as above by a link $K\cup L$, but now in $M\times\{\varepsilon\}$,
so we write $X = M[K,L]$.  
This presents $\pi_1(X)$ as a quotient of the free product $\pi_1(M)*F(x_1,\dots,x_k)$ (again the $x_i$ are meridians of $K$) by the normal subgroup generated by the elements $r_1,\dots, r_\ell$ (determined up to conjugacy) traced out by the components of $L$.  To indicate this, we write   
$$
\pi_1(X) \ = \ (\pi_1(M),  \, x_1, \, \dots, \, x_k \st r_1, \, \dots, \, r_\ell)\,.
$$
 
\begin{definition}\label{accobordism}
A relative handle structure $X=M[K,L]$ as above is called an \bit{\ac\,structure} if it satisfies the conditions of \defref{acmanifold}, that is, $r_i=x_i$ for $i\le n$ and $r_i=1$ for $i>n$ for some $n \leq \min(k,\ell)$ 
and a suitable ordering of the handles.  The quotient $X/M$ (crushing $M$ to a point) is then homotopy equivalent to a wedge of circles and $2$-spheres: just circles for type \ione\ structures (where $n=\ell$) and just $2$-spheres for those of type \itwo\ ($n=k$).
Note that for the latter type, the map $\pi_1(M) \to \pi_1(X)$ induced by inclusion is an isomorphism.

The cobordism $X$  is  \bit{\ac}, or an \bit{\ac\,cobordism} from $M$ to $N$, if it
admits an \ac\,structure.  
In particular, if $X$ can be built from $M\times [0,1]$ by adding an equal number of homotopically cancelling $1$ and $2$-handles, and so is both of type \ione\ and \itwo, then we call it an \bit{\ac\ homology cobordism}.  
\end{definition}

\begin{remark}\label{pi1}
The ``homology cobordism" condition means that the inclusions of both $M$ and $N$ into $X$ induce isomorphisms on homology, so in particular $M$ is a homology sphere if and only if $N$ is. This follows from the \ac\,condition using the homology exact sequences of $(X,M)$ and $(X,N)$, since the contractibility of $X/M$ implies the vanishing of $H_*(X,M)$, and so also $H_*(X,N)$ by duality.    
\end{remark}

\headtop{Encasement and tight embeddings}

Embeddings of one \ac\,manifold in another will frequently arise in our subsequent discussions. It is often helpful when  \ac\,structures can be chosen compatibly, and also when the complement of the embedded manifold has a suitable \ac\,structure. 
For these purposes, we formulate the notions of ``encasing" one \ac\,manifold in another, and of ``tight" embeddings.   

\begin{definition}\label{encase} 
An \ac\,manifold $X$ embedded in the interior of an \ac\,manifold $Y$ is said to be {\it encased} in $Y$ if there is an \ac\,structure on $X$ 
that extends to one on $Y$.  
\end{definition}

\begin{remark}\label{twistAC}
When $X$ is a cork encased in $Y$, then in fact {\it every} \ac\,structure on $X$ will extend to $Y$ and {\it every}  cork twist $Y_{X,h}$ will still be \ac.  This can be seen since an \ac\,structure on $Y$ extending one on $X$ induces presentations $\pi_1 Y = (x_1,\dots ,x_k \st$ $r_1,\dots r_\ell)$ and $\pi_1 X = (x_1,\dots ,x_m \st$ $r_1,\dots r_m)$ where $r_i = x_i$ if $i \leq n$ and $r_i=1$ if $i>n$, for some $n \leq \min(k,\ell)$ and $m \leq n$. Changing or twisting the \ac\,structure on $X$ while fixing the handles of $Y-X$ gives a new presentation for $\pi_1 Y$ and $\pi_1(Y_{X,h})$, respectively, where the relations coming from the $2$-handles of $Y-X$ are changed only up to generators coming from the $1$-handles of $X$. However, by assumption, the $2$-handles of $X$ still homotopically cancel the $1$-handles of $X$. So, the $2$-handles of $Y-X$ may be slid over those in $X$ to be attached along their original relations, giving an \ac\,structure  on $Y$. 
\end{remark}

\begin{definition}\label{tight}
Let $X$ and $Y$ be compact $4$-manifolds, and assume that $Y$ is simply-connected with connected (possibly empty) boundary. A \bit{tight embedding} of $X$ in $Y$ will mean an embedding $X\subset\interior(Y)$ for which $Y-\interior(X)$ is either an \ac\ manifold when $Y$ is closed, or an \ac\,cobordism from $\del Y$ to $\del X$ when $Y$ has boundary. When we say ``$X\subset Y$ is tight", it is understood that $Y$ is simply-connected and the embedding is in $\interior(Y)$, and if $X$ is a cork (with or without a diffeomorphism on its boundary) we may just refer to it as a ``tight cork in $Y$".   
\end{definition}       

\begin{remark}\label{1connected}  If $X$ is a tight cork in $Y$, then the complement $Y-\interior(X)$ has trivial first homology, so must be \ac\ of type \itwo, i.e.\ built from either a $4$-ball or a collar on $\del Y$ with only $1$ and $2$-handles so that some subset of the $2$-handles homotopically cancel the $1$-handles. Furthermore, $\pi_1(Y-X)$ will be trivial when $Y$ is closed.  If $Y$ has boundary, then the inclusion $\del Y \hto Y-X$ is surjective on $\pi_1$, as noted in \defref{accobordism}, so if $Y$ is itself embedded with simply-connected complement in a closed $4$-manifold $Z$, then $Z-X$ will again be simply-connected  (by Van Kampen's Theorem). 
\end{remark}

\headtop{Multicorks}

For the present purposes, it is useful to broaden the notion of corks and their cork twists.

\begin{definition}\label{multicork}
A \bit{multicork} is a list $\calc = (C_0,C_1,\dots)$ of compact contractible $4$-manifolds equipped with boundary diffeomorphisms $h_i\colon\del C_i\to \del C_0$, with $h_0 = \id$.  The $C_i$ are called the \bit{components} of $\calc$.   The \bit{order} of $\calc$ is the length of the list, which may be finite or infinite.  The boundaries of the $C_i$ are by definition all diffeomorphic, but the $C_i$ themselves need not be  \cite{akbulut-ruberman:absolute} (although they must be homeomorphic \cite{freedman:simply-connected}; see the footnote to \ref{ICT} below).  If the $C_i$ are all \ac, then $\calc$ is called an \bit{\ac\,multicork}.  If in addition each ``twisted double" $-C_0\cup_{h_i}C_i$ is diffeomorphic to the $4$-sphere, then $\calc$ is called a \bit{simple} multicork.
\end{definition}

{\it Multicorks generalize corks.}  Indeed, there is an order-preserving (non-surjective) injection 
$$
\mu\colon\text{Corks} \lhto \text{Multicorks}
$$
sending any cork $(C,h)$ to the constant multicork $\mu(C,h) = (C,C,\dots)$ of the same order with boundary identifications $h_{i} = h^{i}$.  Note that if $(C,h)$ is simple, then so is $\mu(C,h)$.

\begin{definition}\label{corkreplacements}
An \bit{embedding} of a multicork $\calc = (C_0,C_1,\dots)$ in a 4-manifold $X$ is an embedding $C_0\subset\interior(X)$, and if this embedding is tight, we simply say ``$\calc$ is a tight multicork in $X$".  Such an embedding generates a list of \bit{cork replacements} 
$$
X_{C_0,C_i} \ := \ (X-\interior(C_0))\ \cup_{h_{i}}\,C_i.
$$
These replacements depend of course on the embedding $C_0\subset X$ and the boundary identifications $h_i$, but this dependence is suppressed in the notation since it can often be gleaned from the context.  Note that $\del X_{C_0,C_i}$ and $\del X$ are naturally identified, and will be viewed as being literally equal.  
\end{definition}

{\it Cork replacements generalize cork twists.}  Indeed, the cork twists $X_{C,h^i}$ of a cork $(C,h)$ associated with an embedding $C\subset X$ are the cork replacements $X_{C_0,C_{i}}$ of the multicork $ \mu(C,h)$ associated with the same embedding.  It will be important for what follows to know when we can go in the opposite direction: Given a multicork $\calc$ in $X$, under what conditions is there a cork in $X$ of the same order whose twists give manifolds diffeomorphic to the cork replacements of $\calc$?

\begin{definition}\label{correlated} An embedded cork $(C,h)$ in a $4$-manifold $X$ is \bit{correlated} with an embedded multicork $\calc = (C_0,C_1,\dots)$ in $X$ if the cork twists $X_{C,h^i}$ are diffeomorphic rel boundary to the cork replacements $X_{C_0,C_i}$ for all $i$ (recall that the boundaries of $X_{C,h^i}$ and $X_{C_0,C_i}$ are naturally identified with $\del X$, and thus with each other).  Similarly $(C,h)$ is \bit{correlated} with a list $(A_i,\tau_i)$ of involutory corks embedded in $X$ if the $X_{C,h^i}$ are diffeomorphic rel boundary to $X_{A_i,\tau_i}$ for each $i$.
\end{definition}

The ``pinwheel lemma" below shows that any finite simple multicork $\calc$ has an associated higher order cork which has embeddings correlated with any embedding of $\calc$.  This will be used in \S2 to prove a related result -- the ``consolidation theorem" -- for finite lists of simple involutory corks.

\head{Pinwheels}

\def\rot{\pi}

\begin{definition}\label{def:pinwheel}
The \bit{pinwheel} of a finite multicork $\calc = (C_0,\dots,C_{n-1})$ is the finite cork 
$(P,\rot)$,
where $P$ is the boundary connected sum 
$
C_0\bcs\cdots\bcs C_{n-1}
$
and $\rot$ is the ``obvious" boundary rotation shifting each $\del C_i$ to $\del C_{i-1}$ (here subscripts are taken mod $n$).  More precisely, fix a linear rotation $\rot$ of the 3-sphere of order $n$ and a principal orbit $x_0, \dots,x_{n-1}$ of the action of the cyclic group generated by $\rot$, with $x_i = \rot^{-i}(x_0)$ for each $i$. Choose $y_0\in\del C_0$, and set $y_i = h_i^{-1}(y_0)\in \del C_i$, where $h_{i}\colon\del C_i\to \del C_0$ are the boundary identifications of $\calc$.  Now build $\bsum$ from the disjoint union $\ddsum$ by adding 1-handles joining $x_i\in S^3$ to $y_i\in \del C_i$ for each $i$.  Then $\rot$ extends to a rotation of 
$$
\del P \ = \ \csum \ \cong \ \del C_0\cs\cdots\cs\del C_0
$$
of order $n$, sending each $\del C_i$ to $\del C_{i-1}$ via $h_{i-1}^{-1}h_i$, as shown in \figref{pinwheels}a for the case $n=4$. 
\end{definition}

\begin{remark}\label{pwprop} 
The pinwheel construction preserves many of the properties discussed above. For example, it follows from the definitions that the pinwheel of an \ac\ multicork is an \ac\ cork, and the pinwheel of a simple multicork of order $2$ must also be simple. It is unclear whether or not this last fact still holds for multicorks of higher order.  
\end{remark}

\fig{145}{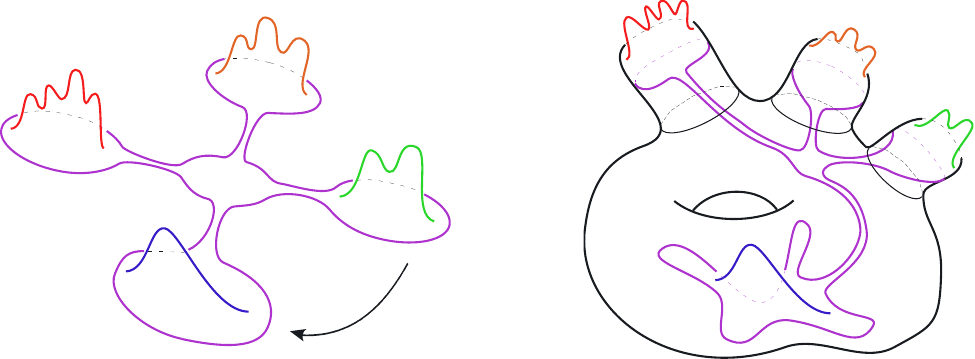}{
\put(-310,68){$B^4$}
\put(-330,22){$C_0$}
\put(-243,57){$C_1$}
\put(-294,107){$C_2$}
\put(-375,85){$C_3$}
\put(-238,15){$\rot$}
\put(-147,100){$B_3$}
\put(-88,120){$B_2$}
\put(-10,65){$B_1$}
\put(-380,-15){a) Pinwheel of an order $4$ multicork}
\put(-150,-15){b) Embedding the pinwheel}
\caption{Pinwheels}
\label{pinwheels}}

Given an embedding of a finite simple multicork $\calc$, the ``Pinwheel Lemma" below produces an embedding of the pinwheel of $\calc$ 
whose cork twists are diffeomorphic (rel boundary) to the cork replacements of $\calc$, that is, the embeddings of $\calc$ and 
its pinwheel
are correlated (see \defref{correlated}). 
The proof of this lemma generalizes an argument of Matveyev \cite{matveyev} (see also Kirby \cite{kirby:cork}) that produces involutory cork embeddings correlated with embeddings of order 2 multicorks; a similar generalization was used in the construction of equivariant corks in \cite{akmr:equivariant}.


\begin{namedtheorem}{Pinwheel Lemma}\label{pinwheel}
For any embedding of a finite simple multicork $\,\calc$ in a $4$-manifold $X$, there is a correlated embedding of its pinwheel $(P,\rot)$ in $X$.  If the embedding of $\calc$ is tight, then the embedding of $P$ can also be chosen to be tight.  
\end{namedtheorem}   


\proof  
Let $\calc = (C_0,\dots,C_{n-1})$ with $C_0\subset X$ and boundary diffeomorphisms $h_i\colon\del C_i\to\del C_0$. 
Since $\calc$ is simple, each $-C_0 \cup_{h_i}C_i\cong S^4$.  Thus the $C_i$ for $i>0$ embed in disjoint 4-balls $B_i$ in $X-C_0$ with $B_i-C_i$ diffeomorphic rel boundary to a punctured copy of $\interior(-C_0)$.  These embeddings are tight since $C_0$ is \ac. 
Combining them with the original embedding $C_0\subset X$ gives an embedding of $\Q = C_0\sqcup\cdots\sqcup C_{n-1}$ in $X$, which then extends 
to an embedding of $P=\bsum$ in $X$, guided by a collection of disjoint arcs in $X-\Q$ as illustrated in \figref{pinwheels}b. 

For each $i$, the map $\rot^i$ rotates $P$ (clockwise in the figure) by a $(i/n){\text{th}}$ turn, and hence the cork twist $X_{P,\rot^i}$ is diffeomorphic rel boundary to the 4-manifold $X_i$ obtained from $X$ by replacing each $C_j$ by $C_{i+j}$ (with subscripts mod $n$).  These replacements have no effect except when $j=0$ since $\calc$ is simple (and all diffeomorphisms of $\del B^4$ extend to $B^4$ by Cerf \cite{cerf}) so $X_i$ is diffeomorphic rel boundary to $X_{C_0,C_i}$.  Therefore the embedding $P\subset X$ is correlated with the embedding of $\calc$ in $X$. 

Now assume $C_0\subset X$ is tight.  Recall that this presumes that $X$ is compact and simply-connected with connected, possibly empty, boundary.  In fact we can assume that $\del X$ is nonempty since the closed case follows by removing a $4$-ball.  Then, the embedding of $C = \esum$ in $X$ is tight, since $C$ is an isotopic deformation of $C_0$.  But $X-\interior(P)$ is obtained from $X-\interior(C)$ by attaching only homotopically cancelling $1$ and $2$-handles, since each $B_i-\interior(C_i)$
is an \ac\ homology cobordism.  Thus $X-\interior(P)$ is an \ac\,cobordism from $\del X$ to $\del P$. That is, $P\subset X$ is tight.
\qed 

\vskip-5pt
\vskip-5pt


\head{The Involutory Cork Theorem}

The  Involutory Cork Theorem, due to Curtis-Freedman-Hsiang-Stong \cite{curtis-freedman-hsiang-stong} and Matveyev \cite{matveyev}, states that any two (smooth) closed simply-connected $4$-manifolds that are homotopy equivalent, and thus homeomorphic, are related by a single  cork replacement (see \defref{corkreplacements}).   It was shown in \cite{curtis-freedman-hsiang-stong}, based on earlier work of Stong \cite{stong}, that the common complement of the corks can be chosen to be simply-connected.    Moreover, as demonstrated in \cite{matveyev}, this cork replacement can be accomplished by an involutory cork twist. A nice exposition of the full theorem is given in \cite{kirby:cork} 

We now state a relative version of the theorem in a form that will be convenient for our purposes; it immediately implies the absolute version by removing a $4$-ball.  The proof can be gleaned from a careful reading of \cite{curtis-freedman-hsiang-stong} and \cite{kirby:cork}, but for the reader's convenience we include a sketch.  

\begin{namedtheorem}{Relative Involutory Cork Theorem} \label{ICT}  
If $X_0$ and $X_1$ are homotopy equivalent, smooth, compact simply-connected $4\text{-manifolds}$ with diffeomorphic homology sphere boundaries, then any diffeomorphism $f\colon \del X_0 \to \del X_1$ extends to a diffeomorphism $(X_0)_{P,\tau} \to X_1$ for some tight simple involutory cork $(P,\tau)$ in $X_0$.\foot{In particular $f$ extends to a homeomorphism $X_0\to X_1$ (known previously to Freedman \cite{freedman:simply-connected}; see also Boyer \cite{boyer}).} \end{namedtheorem}

\vskip-10pt
\vskip-10pt

\proof First construct a relative h-cobordism $W$ from $X_0$ to $X_1$ with the mapping cylinder $M_f$ on the lateral boundary, as follows: The closed manifold $X_1 \cup_f -X_0$ has zero signature, hence bounds a 5-manifold which, appropriately surgered, has a relative handlebody structure with only $2$ and $3\text{-handles}$, as in Wall's foundational paper \cite{wall:4-manifolds}.   Splitting this handlebody along its middle level (between the 2 and $3\text{-handles}$) and re-gluing using Wall's Theorem 2 \cite[pg.\ 136]{wall:diffeomorphisms} (note that it may be necessary to introduce an extra cancelling 2/3-handle pair to apply Wall's result) gives the desired h-cobordism $W$, built from $X_0\times [0,\varepsilon]$ by adding only $2$ and $3$-handles.  

The middle level $X_{1/2}$ of the $h$-cobordism $W$ contains two sets of embedded 2-spheres, the belt spheres $S_-$ of the 2-handles and the attaching spheres $S_+$ of the 3-handles of $W$.  These spheres can be assumed to meet transversely, to pair algebraically, and after a suitable sequence of finger moves as pioneered by Casson \cite{casson} in the 1970s, to have connected union with simply-connected complement $X_{1/2}-(S_-\cup S_+)$.  

Let $U_{1/2}$ be a regular neighborhood of $S_-\cup S_+ \cup T$ in $X_{1/2}$, where $T$ is an embedded tree joining generic base points in the spheres in $S_-\cup S_+$ to a base point in their complement.
Surgering $U_{1/2}$ along $S_-$ and $S_+$ in turn 
yields submanifolds $U_0\subset X_0$ and $U_1\subset X_1$ joined by a cobordism $U\subset W$ built from a thickening $U_{1/2}\times [a,b]$ 
of the middle level by adding 
$3$-handles along the spheres $S = (S_- \times \{a\}) \cup (S_+\times\{b\})$.  Working up from the bottom, this gives a relative handle structure for $U$ (adding $2$ and $3$-handles to a thickening of $U_0$) that extends trivially to one on $W$. 

Choose a Morse function $h\colon W\to[0,1]$ with $X_i = h^{-1}(i)$ for $i=0$, $1/2$ and $1$ that is compatible with this handle structure on $W$ and that extends the standard height function on $M_f$ (see \figref{ict}).  Note that $W-U$ acquires a product structure from the gradient flow of $h$ (with respect to a suitable metric) which gives a diffeomorphism 
$
G\colon X_0\!-U_0 \to X_1\!-U_1
$ 
extending $f$ on the boundary. 

\fig{150}{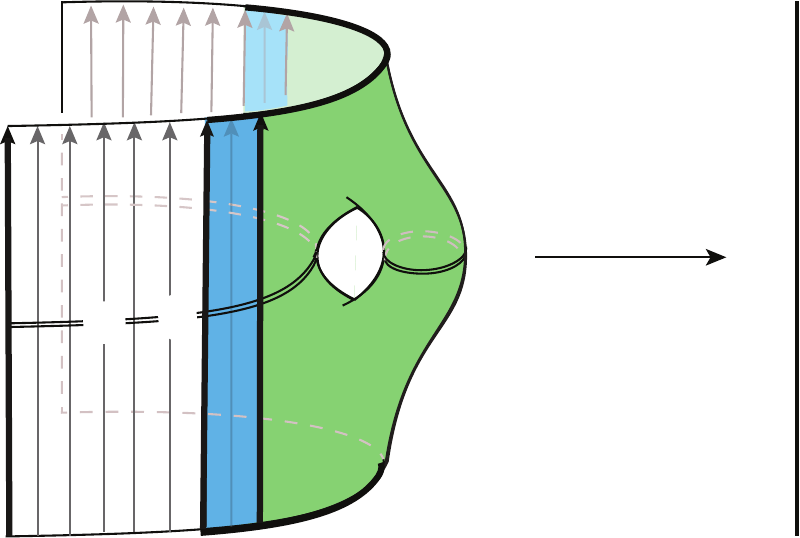}{
\put(-240,20){\large$W$}
\put(-130,50){\large$U$}
\put(-163,40) {\large$V$}
\put(-110,25) {\large$C = U\cup V$}
\put(-230,56){$f$}
\put(-175,59){$g$}
\put(-198,57){$G$}
\put(-140,-5){\small$C_0$}
\put(-139,150){\small$C_1$}
\put(-52,65){$h$}
\put(4,1){$0$}
\put(4,142){$1$}
\put(4,75){$1/2$}
\caption{Proof of the Relative Involutory Cork Theorem}
\label{ict}}

Now the handlebody techniques in \cite{stong} can be used to encase $U_{1/2}$ (as in Definition \ref{encase}) in a tight \ac\,submanifold $C_{1/2}$ of $X_{1/2}$ of type \itwo\ with $H_2(C_{1/2}) \cong H_2(U_{1/2})$ (both groups are generated by the spheres in $S_-$ and $S_+$); this is immediate from the Encasement \lemref{encasement} in the next section (see \remref{ICTproof}).    Let $V$ be the product cobordism swept out by $C_{1/2} - \interior(U_{1/2})$ under the gradient flow of $h$, and set $C = U\cup V$.  For $i=0,1$, the manifold $C_i = C\cap X_i$ is an \ac\,cork in $X_i$ containing $U_i$. Although the manifolds $U_0$ and $U_1$ are not necessarily \ac\, (see \cite[Figure 2]{kirby:cork} for example), it follows from the proof of the Encasement Lemma applied to $U_{1/2}$ that all $1$-handles of $U_0$ and $U_1$ are indeed algebraically cancelled by $2$-handles in $C_{1/2}$. Thus the $C_i$ are \ac\, corks, both tightly embedded since $C_{1/2}$ is. Their boundaries $\del C_0$ and $\del C_1$ are identified via the restriction $g=G|_{\del C_0}$, as indicated in Figure \ref{ict}.  This makes $(C_0,C_1)$ into an \ac\,multicork in $X_0$, such that the restriction of $G$ to $X_0\!-\!C_0$ extends to a diffeomorphism $(X_0)_{C_0,C_1} \to X_1$ extending $f$ on the boundary.

In fact the multicork $(C_0,C_1)$ is {\it simple} (see \defref{multicork}), i.e.\ both $C_0\cup_{\id}-C_0$ and $\Sigma = C_1\cup_{g}-C_0$ are diffeomorphic to a $4$-sphere.  Indeed the double of $C_0$ is standard by Remark \ref{andrewscurtis} \ib, while $\Sigma = \del C$ and $C$ is easily seen to be a $5$-ball (cf.\ Kirby \cite{kirby:cork}):  It is built from $C_{1/2} \times [a,b]$ by adding $3$-handles along the spheres $S = (S_- \times \{a\}) \cup (S_+\times\{b\})$.  Since $C_{1/2}$ is the result of surgering a collection of curves in an \ac\ cork (either $C_0$ or $C_1$), its thickening $C_{1/2} \times [a,b]$ can be built from the $5$-ball by adding only $2$-handles.  But the core of each $2$-handle lies in one of the spheres in $S$, and so the $2$-handle is cancelled by the $3$-handle of $C$ attached to that sphere.  

Since the multicork $(C_0,C_1)$ is simple and tightly embedded, its involutory pinwheel $(P,\rot)$ is also simple 
and has a tight embedding $P \subset X_0$  correlated with the multicork embedding $C_0\subset X_0$, by \remref{pwprop} and \lemref{pinwheel}. Hence $f$ extends to a diffeomorphism $(X_0)_{P, \tau} \to X_1$ as desired.   
\qed

\begin{remark*}
Theorem \ref{ICT} need not hold when $\del X_0$ is not a homology sphere. For example, setting $X = X_0 = X_1 = S^2\!\times\!D^2$, the Gluck twist $\tau$ on $\del X$ from \cite{gluck:2-spheres} does not extend across $X$, even after a cork twist, as $X\cup_\id -X = \sss$ and $X\cup_\tau -X = \sts$ are not even homotopy equivalent. 
\end{remark*}


\begin{corollary}\label{twisteddouble}
Every homotopy $4$-sphere $\Sigma$ is diffeomorphic to a twisted double $W \cup_{h} -W$ of some \ac\,cork $(W,h)$. 

\end{corollary}

\proof Decompose $\Sigma$ as the union of a homotopy 4-ball $X$ and a 4-ball $B$ meeting in their common 3-sphere boundary $S$.  By \thmref{ICT}, the identity map $\del X \to \del B$ extends to a diffeomorphism $f\colon X_{P,\tau}\to B$ for some simple involutory cork $(P, \tau)$ tightly embedded in $X$.  Enlarging the corks $P\subset X$ and $Q = f(P)\subset B$ along arcs to the boundary, we may assume that $P$ and $Q$ (shown in green on the left side of \figref{double}) meet $S$ in disjoint 3-balls.  Set $P'=\cl(X-P)$ and $Q'=\cl(B-Q)$ (shown in purple and blue on the left side).

\fig{120}{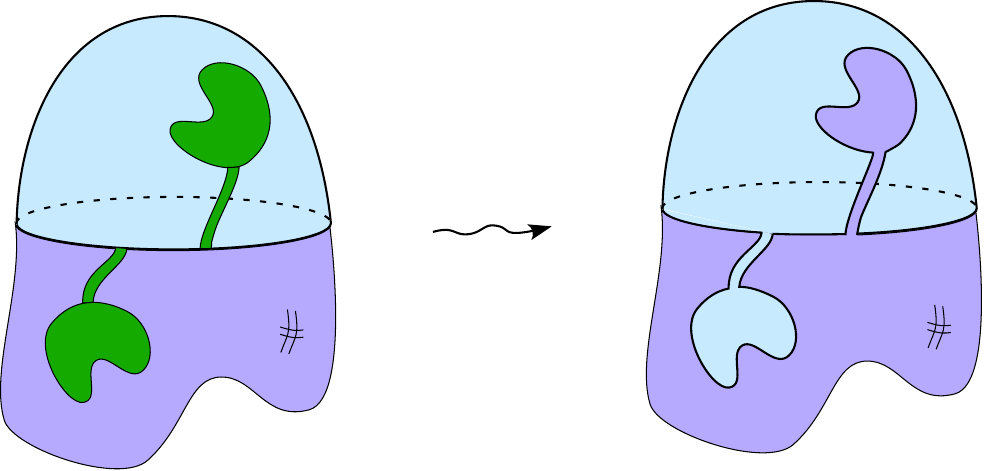}{
\put(-253,90){$B$}
\put(-263,26){$X$}
\put(-258,60){$S$}
\put(-230,32){$P$}
\put(-195,85){\small$Q$}
\put(2,34){$W$}
\put(-4,90){$Z\cong -W$}
\caption{Both decompositions of $\Sigma$}
\label{double}}

\noindent Then $P'$ and $Q'$ are diffeomorphic \ac\,corks, each built from a $4$-ball (the collar on a punctured copy of $S$) by adding homotopically cancelling $1$ and $2$-handles.   Hence $\Sigma = X\cup B = W \cup Z$, where $W = P' \bcs Q$ and $Z = Q' \bcs P$ (shown in purple and blue on the right side of the figure) are \ac\,corks.  But $f$ induces a diffeomorphism $g\colon\!\!-W\to Z$, and so $\Sigma = W\cup_h-W$ where $h=g\vert_{-\partial W}$. \qed

\begin{remark*}
This result extends to any homeomorphic pair of closed simply-connected $4$-manifolds $X$ and $Y$, showing that $X \# -Y$ is diffeomorphic to a twisted double $W \cup_h -W$, where the honest double $W\cup -W$ is diffeomorphic to a connected sum of $S^2$-bundles over $S^2$. 
\end{remark*}

\section{Consolidating Corks}\label{fcs}


The main tool for proving the cork theorems in the introduction is the following result, showing that in suitable $4$-manifolds $X$, any finite collection of simple involutory corks in $X$ is correlated with a single higher order cork in $X$:

\begin{namedtheorem}{Consolidation Theorem}\label{consolidation} 
Let $(A_1,\tau_1),\dots,(A_{n-1},\tau_{n-1})$ be simple involutory corks embedded in a compact simply-connected $4$-manifold $X$ whose boundary is either empty or a homology sphere.  Then the $A_i$ can be isotoped so that they lie in a single \ac\ cork $(C,h)$ of order $n$ in $X$ whose twists $X_{C,h^i}$ are diffeomorphic rel boundary to $X_{A_i,\tau_i}$ for each $i$.  If $X$ is closed and the $A_i$ have simply-connected complements, then $C\subset X$ may be chosen with simply-connected complement. 
\end{namedtheorem}  

\noindent The proof of this theorem requires careful manipulation of a handle structure for $X$.  We begin by stating some generalities about such structures, and establish two preliminary results: the ``encasement" and ``finger" lemmas.

\head{Handle Manipulations}

Consider the case when $\del X$ is a homology sphere, or more generally nonempty and connected, and fix a handle structure on $X$ with a single $0$-handle and no $4$-handles.  By convention, assume that no collars are added between the handles, except for a final collar after all of the handles are attached.  The \bit{$k$-skeleton} $\ks$ of $X$ consists of all the handles of index $\le k$.  Inverting $X$ produces the dual handlebody $\X$ with \bit{dual $k$-skeleton} $\dualks$, consisting of a collar on $\del X$ with all the dual handles of index $\le k$ (upside down handles of index $\ge 4-k$) attached.


The \bit{middle level} of $X$ is the bounded 3-manifold $M \ = \ \del\one\cap\del\dualone$. 
Because of our collar convention, $M$ can be seen as the complement in $\del\one$ of the $2$-handle attaching regions, or equivalently, the complement in $\del\dualone-\del X$ of the dual $2$-handle attaching regions.  By general position, it follows that the inclusion induced maps 
$$
\pi_1(\one) \underset p\lfrom \pi_1(M) \underset \dualp\lto \pi_1(\dualone)
$$
are surjective, where a base point in $M$ is implicit.  Note that $\pi_1(\one)$ is the free group $F(x_1,\dots,x_k)$, where the $x_i$ are given by the $1$-handles of $X$, and $\pi_1(\dualone)$ is the free product of $\pi_1(\del X)$ with the free group $F(x^*_1,\dots,\,x^*_m)$ on generators $x^*_i$ given by the dual 1-handles of $X$.  

The meridians and framed longitudes of the attaching circles of the 2-handles $H_i$ of $X$ are unbased loops in $M$, and thus correspond to conjugacy classes $m_i$ and $\ell_i$ in $\pi_1(M)$.  It is geometrically evident that the projected conjugacy classes $p(m_i)$ and $\dualp(\ell_i)$ are trivial, while the classes
$$
r_i \ = \ p(\ell_i)
\quad\text{and}\quad
\dri \ = \ \dualp(m_i)
$$
recording the attaching circles of the $H_i$ and their duals $\H_i$ are typically nontrivial.  These classes (or their representatives) are called the \bit{relators} and \bit{dual relators} of $X$.  Since $\pi_1(X)$ is trivial, they normally generate $\pi_1(\one)$ and $\pi_1(\dualone)$, respectively, and it follows that the map
$$
p \times \dualp \,:\,  \pi_1(M) \lto \pi_1(\one) \times \pi_1(\dualone) \ = \ F(x_1,\dots,x_k) \times (\pi_1(\del X)*F(x^*_1,\dots,\,x^*_m))
$$
is surjective; cf.\ Stong's Lemma \cite[p.\,500]{stong} in the closed case (the proof is the same here).  This was a key observation from \cite{stong} used to produce corks with simply-connected complements in closed $4$-manifolds \cite{curtis-freedman-hsiang-stong}, and is essential to our argument below.  

The other input from \cite{stong} is an analysis of the effect on the relators and dual relators of sliding one 2-handle $H_i$ over another $H_j$ along a path $\gamma$ joining their attaching circles \cite[pp.\,499-500]{stong}.  Stong notes that in the dual picture, $\H_j$ slides over $-\H_i$ along the reverse of $\gamma$, as illustrated in \figref{handleslide} by comparing longitudes and meridians of the attaching circles before and after the slide.  

The path $\gamma$ determines an element in $\pi_1(M)$ (also denoted $\gamma$ by abuse of notation) via a choice of base paths from its endpoints to the base point in $M$.  Set $(c,c^*) = (p(\gamma),\dualp(\gamma))$.  Note that the effect of the handleslide on the relators and dual relatiors depends only on $c$ and $c^*$, and so we simply say ``slide $H_i$ over $H_j$ along $(c,c^*)$"  to specify its relevant features.  In particular, the slide changes $r_i$ to $r_icr_j\overbar c$, and thus $\drj$ to $\drj\overbar c^*\odri c^*$ by Stong's observation, fixing all the other relators and dual relators; here and below overbars denote inverses.
This is demonstrated in \figref{handleslide}. Since $p\times\dualp$ is onto, $(c,c^*) \in \pi_1(\one) \times \pi_1(\dualone)$ can be arbitrary by choosing $\gamma$ appropriately. 

\fig{110}{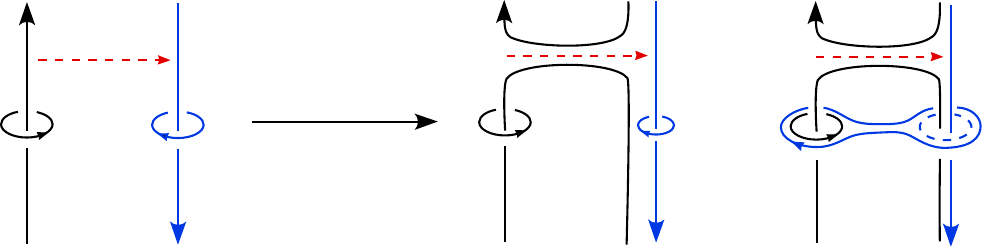}{
\put(-430,-10){$H_i$}
\put(-370,-10){$H_j$}
\put(-446,55){\small$\dri$}
\put(-380,55){\small$\drj$}
\put(-230,-10){$H_i + H_j$}
\put(-150,-10){$H_j$}
\put(-93,-10){$H_i+ H_j$}
\put(-102,15){\small$r_icr_j\overbar c$}
\put(-18,-10){$H_j$}
\put(-120,50){\Large$=$}
\put(-393,90){$\gamma$}
\put(-315,42){\small Handleslide}
\put(-55,35){\small\color{blue}{$\drj\overbar c^*\odri c^*$}}
\caption{Single slide along $\gamma = (c,c^*)$}
\label{handleslide}}

\noindent The following two operations suffice for our purposes:

\medskip
\def\and{\quad\text{and}\quad}

 \noindent {\bf Single slide}\,: {\it Slide $H_i$ over $H_j$ along $(c, 1)$ \,$($or along $(1,c^*))$.  This changes $r_i$ to $r_ic r_j\overbar c$ and $\drj$ to $\drj \odri$ \,$($or $r_i$ to $r_i r_j$ and $\drj$ to $\drj\overbar c^*\odri c^*)$ \,while fixing all other relators and dual relators.}

\noindent {\bf Double slide}\,: {\it Slide $H_i$ over $H_j$ along $(c, 1)$, and then over $-H_j$ along $(1, 1)$.  This changes $r_i$ to $r_i c r_j \overbar c \,\overbar r_j = r_i  [c,r_j]$ while fixing all other relators and dual relators.}

\medskip

Thus one can use {\it single slides} to {\it multiply a relator $($or dual relator$)$ by a conjugate of another}, and {\it double slides} to {\it multiply a relator by a commutator of another}.   It follows from a variation on Stong's argument \cite{stong} (which uses the algebraic observations above to control the fundamental group of the cork complement in the Involutory Cork Theorem \cite{curtis-freedman-hsiang-stong})
that \ac\,submanifolds can often be encased in \ac\,submanifolds of type \itwo\ that are tightly embedded:

\head{Encasement}

\begin{namedtheorem}{Encasement Lemma}\label{encasement} 
Let $X$ be a compact simply-connected $4$-manifold that is either closed or bounded by a homology sphere.  Then any \ac\,manifold $U$ 
lying in the interior of $X$ can be encased $($in the sense of \textup{\defref{encase}}$)$ in an \ac\,submanifold $V$ of type \itwo\ with $H_2(U) \cong H_2(V)$. If $X-U$ is simply-connected, then the embedding $V\subset X$ can be chosen to be tight $($see \textup{\defref{tight}}$)$.
\end{namedtheorem}

\begin{remark}\label{ICTproof}
The proof will show in fact  that {\it every} \ac\,structure on $U$ extends to one of type \itwo\ on some $V$ containing $U$ with $H_2(V) = H_2(U)$.  This lemma generalizes Stong's Structure Theorem \cite[Theorem 2]{stong} (which is the case when $U=B^4$ and $X$ is closed), and is a key tool in the proof of the Relative Involutory Cork \thmref{ICT} where the \ac\, submanifold $U_{1/2}$ plays the role of $U$.
\end{remark}

\vskip-2pt
\vskip-2pt
\pf
The closed case follows from the bounded case by removing a $4$-ball, so we assume that $\del X$ is a homology sphere.  By hypothesis there is an \ac\,handle structure 
$$
U \ = \  [K_\circ \cup K_-,\,L_\circ\cup L_+]
$$
where the homotopically cancelling $1$-handles and $2$-handles are given by $K_\circ$ and $L_\circ$ (the matching subscripts indicate this cancellation) and the remaining $1$-handles and $2$-handles by $K_-$ and $L_+$. The handles in $U$ will not be moved in subsequent constructions.  See the left side of \figref{fig:encasement} for a concrete example, where $K_\circ \cup L_\circ$ is shown in blue and $K_- \cup L_+$ in green.

\fig{100}{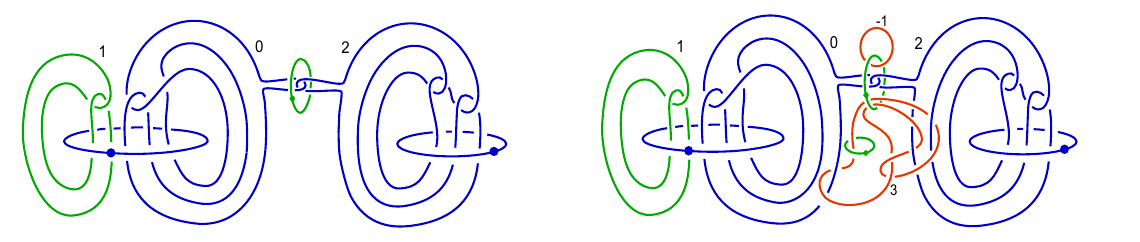}{
\put(-465,10){\small{$L_+$}}
\put(-352,42){\small{$K_-$}}
\put(-365,18){\small{$L_\circ$}}
\put(-338,18){\small{$L_\circ$}}
\put(-395,48){\small{$K_\circ$}}
\put(-308,48){\small{$K_\circ$}}
\caption{$U = [K_\circ \cup K_-,\,L_\circ \cup L_+]$ (left) \ and \ $V = [K_\circ\cup K_1,\,L_\circ\cup L_1\cup L_+]$ (right).}
\label{fig:encasement}}

To construct $V$, extend this handle structure on $U$ to one on all of $X$ with no new $0$ or $4$-handles. This adds a dotted unlink $K$ and a framed link $L$ to the picture, so 
$$
X \ = \ [K_\circ\cup K_-\cup K, L_\circ\cup L_+\cup L]
$$ 
with some $3$-handles attached.  Set $K_1 = K_-\cup K$.  
Since $\pi_1(X)=1$, we can arrange for some sublink $L_1$ of $L$ to homotopically cancel $K_1$, after first sliding the $2$-handles (attached along) $L$ over one another and over $L_0$, and possibly introducing cancelling 2/3-handle pairs along the way (enlarging $L$) to avoid Andrews-Curtis issues.  Now provisionally set
$$
V \ = \ [K_\circ\cup  K_1, \ L_\circ\cup L_1\cup L_+],
$$
consisting of all the $1$-handles in $X$, an equal number of $2$-handles homotopically cancelling those $1$-handles, and some extra $2$-handles (from $U$) attached homotopically trivially (to $L_+$).  Thus $V$ is \ac\ of type \itwo\  encasing $U$ with $H_2(U) \cong H_2(V)$, proving the first statement in the lemma.   See the right side of \figref{fig:encasement}, where $K_\circ\cup L_\circ$ is blue, $K_1 \cup L_+$ is green, and $L_1$ is orange.

Now assume $\pi_1(X-U)=1$.  We will show how to modify $V$ to make it tight, leaving $U$ fixed.  In fact the entire handlebody $[K_\circ      \cup K_1,L_\circ \cup L_+]$ will remain untouched except that 
$K_1$ may expand to include some new $1$-handles, whereas $L$ (and thus $L_1$) may change. 

To begin the argument, first note that the relative handle structure on $X-\interior(U)$, built from $\del U$ by attaching $1$-handles, $2$-handles (along $L$), and $3$-handles, gives a dual handle structure on the same manifold built from $\del X$ by attaching dual $1$-handles (the upside down $3$-handles), dual $2$-handles (along $L^*$, the meridians of $L$), and dual $3$-handles.  This dual structure yields a presentation of the trivial group   
$$
1 \ = \ \pi_1(X-U) \ = \ (\pi_1(\del X),\, x_1^*, \,\dots, \,x_m^* \st r^*_1,\, \dots,\, r^*_n)
$$ 
as in the preamble 
to \defref{accobordism}, where the $x^*_i$ correspond to the dual $1$-handles  and the $\dri$ are the relators in the free product $F := \pi_1(\del X)*F(x^*_1,\dots,x^*_m)$ given by the dual $2$-handles $L^*$.  

Now introduce $m$ cancelling dual $2/3$-handle pairs, one for each dual $1$-handle, and slide these new dual $2$-handles over $L^*$ to homotopically cancel the dual $1$-handles. This can be achieved by single slides as defined above since each dual generator $x_j^*$ is a product of conjugates of the $r_i^*$.
Dually (i.e.\ rightside up) the effect on $X$ is to introduce $m$ cancelling $1/2$-pairs with $2$-handles, called the \bit{upper $2$-handles}, attached along a link $J$, and then to slide $L$ over $J$ so that the dual upper $2$-handles homotopically cancel the dual $1$-handles.   

Next extend $K_1$ to include all the new $1$-handles, and consider the subhandlebody
$$
W \ := \ [K_\circ\cup K_1, L_\circ\cup L_+\cup L] \ \subset \ X.
$$
We claim that $H_1(W)$ is trivial.  Indeed $X-\interior(W)$ is an \ac\ homology cobordism from $\del X$ to $\del W$ (see \defref{accobordism}) obtained from a collar on $\del X$ by adding all the dual $1$-handles and their homotopically cancelling dual $2$-handles, attached to the meridians $J^*$ of $J$.  Since $\del X$ is a homology sphere, so is $\del W$ by \remref{pi1},  whence $H_1(W) = 0$ by a Mayer-Vietoris argument.

Observe also that the fundamental group of the subhandlebody $[K_\circ\cup K_1,L_\circ\cup L_+] \subset W$ 
is free on the generators $x_1, \dots, x_k$ given by $K_1$, since $L_0$ homotopically cancels $K_\circ$.  
Since $H_1(W)=0$, the $2$-handles $L$ can be slid among themselves so that some subset \emph{homologically} pairs with $K_1$. The handle structure after these handleslides yields a presentation of the trivial group   
$$
\pi_1(X) \ = \ (x_1, \,\dots, \,x_k \st r_1, \, \dots, \, r_n, \, s_1,\, \dots,\, s_m) \ = \ 1
$$ 
where the $x_i$ are given by the $1$-handles $K_1$, the $r_i$ and $s_i$ are the relators given by $L$ and $J$, respectively, and $r_i= x_i \prod_{j=1}^{\ell_i} [a_{ij},b_{ij}]$ for $i=1,\dots,k$ and suitable $a_{ij}, b_{ij} \in F(x_1, \dots, x_k)$. 

Finally introduce $\ell_1+\cdots+\ell_k$ cancelling $2/3$-handle pairs indexed by pairs $(i,j)$ for $i=1,\dots,k$ and $j=1,\dots,\ell_i$, called the \bit{extra handles}.  Each $b_{ij}$ is a product of conjugates of the relators $r_1, \, \dots, \, r_n, \, s_1,\, \dots,\, s_m$ given by the $2$-handles $L\cup J$, so we can single slide the extra $2$-handles over the $2$-handles $L\cup J$ until the $i\!j^{th}$ one $L_{ij}$ has attaching word $b_{ij}$. Then for each $i=1,\dots k$, double slide the $i^{th}$ $2$-handle in $L$ over $L_{ij}$ along $(a_{ij},1)$ successively for $j = 1,\dots,\ell_i$ until the corresponding relator is $r_i=x_i$.  The result is a presentation
$$
\pi_1(X) \ = \ (x_1, \,\dots, \,x_k \st x_1, \, \dots, \, x_k, \, r_{k+1} \, \dots, \, r_n, \, s_1,\, \dots,\, s_m, \, b_{11}, \dots, b_{k\ell_k}) \ = \ 1
$$ 
where the first $k$ relators correspond to a sublink $L_1$ of $L$ that homotopically cancels $K_1$.  Note that the property that each dual extra $1$-handle is homotopically cancelled by its dual extra $2$-handle is maintained while performing these single and double slides.  Now set 
$$
C \ = \ [K_\circ\cup K_1, L_\circ \cup L_1\cup L_+]
$$
which is the union of $U$ with all of the $1$-handles in $X$ and their homotopically cancelling $2$-handles. The embedding $V \subset X$ is tight since every dual $1$-handle in $X-C$ is homotopically cancelled by a dual $2$-handle. \qed 

\vskip-2pt
\vskip-2pt

\head{Repositioning Corks}

The Encasement Lemma can be used to prove the Consolidation Theorem after repositioning any finite collection of corks with simply-connected complement so that their union also has simply-connected complement.  The precise result we need is the following, proved by a standard argument using finger moves:

\begin{namedtheorem}{Finger Lemma}\label{fingermove} 
\it Any finite collection of \ac\,corks $A_1,\dots,A_n$ embedded in a $4$-manifold $X$ can be repositioned by isotopies so that their union $\q$ is an \ac\,manifold of type \ione\ encasing each of the $A_i$, and so that $X-\q$ is simply-connected provided each $X-A_i$ is simply-connected.
\end{namedtheorem}

\begin{proof} 

Fix \ac\,structures $A_i = [K_i,\,L_i]$ for $i=1,\dots,n$. Pull the $A_i$ apart so that their 1-skeleta lie in disjoint 4-balls and their core 2-skeleta intersect transversely.  Then join the 0-handles of the $A_i$ by 1-handles to a disjoint $4$-ball to form a single 0-handle 
whose boundary contains the links $K_i\cup L_i$, embedded in disjoint 3-balls.  In this position, the union $\q$ is obtained from the boundary sum $\b$ by plumbing the $2$-handles together according to the intersection of their cores, each intersection point having the effect of clasping the relevant components of $\Li$ through a dotted circle. This transforms $\Li$ into a more complicated framed link $L$ with clasps surrounded by dotted circles forming an unlink $J$ disjoint from $K=\Ki$ and homotopically unlinked from $L$.  

The result is an \ac\,handle structure $\q = [J\cup K, \,L]$ of type \ione\ extending the one on each $A_i$, built in a canonical way from $\b$. This is illustrated in \figref{corkunion} for the case $n=2$ where $A_1$ (in green) is the Mazur manifold and $A_2$ (in blue) is the Akbulut-Matveyev ``positron" \cite{akbulut-yasui:corks-plugs}, embedded so that their 2-handles are plumbed three times geometrically, and once algebraically, and the components of $J$ (in red) are the small dotted circles around the clasps.

\fig{80}{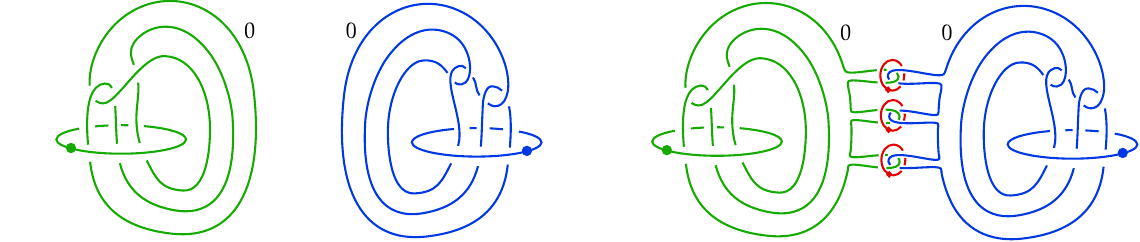}{
\put(-327,-15){$A_1$}
\put(-241,-15){$A_2$}
\put(-130,-15){$A_1\cup A_2  \ \subset \ X$}
\caption{Union of corks}
\label{corkunion}}

To arrange for $\q$ to have simply-connected complement in $X$ when the $A_i$ do, the belt circle of each $2$-handle (say in $A_j$) should bound an immersed disk in $X-\q$.   Since $\pi_1(X-A_j)$ is trivial, such a disk $\Sigma$ can be found in $X-A_j$, but $\Sigma$ might (transversely) intersect the core disk $D$ of some $2$-handle $H$ in another $A_k$, as shown in \figref{finger}a; the (red) dot in the middle represents the $1$-handle $h$ homotopically cancelled by $H$, while all the edges in the figure are cores of $2$-handles. In fact, $\Sigma$ can be assumed to have zero {\sl algebraic} intersection number with $D$, since any point in $\Sigma\cap D$ can be exchanged for some number of cancelling pairs of intersection points between $\Sigma$ and the $2$-handle cores of $A_k$, by a finger move of $\Sigma$ along $D$ across $h$ (see \figref{finger}b).  

\fig{95}{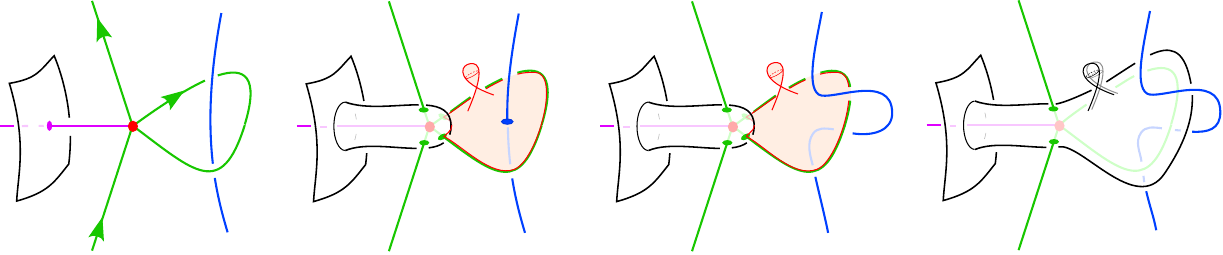}{
\put(-442,-20){a) $\Sigma\,\subset\, X\!-\!A_j$ \hskip 20pt b) finger move along $D$ \hskip 15pt c) finger move along $\Delta$ \hskip 15pt d) Whitney move}
\put(-449,10){\small $\Sigma$}
\put(-426,52){\small $D$}
\put(-409,37){\small $h$}
\put(-338,10){\small $\Sigma$}
\put(-280,40){\small $\Delta$}
\put(-292,58){\small $p$}
\put(-292,35){\small $q$}
\put(-255,33){\small $E$}
\put(-260,80){\small $F$}
\put(-225,10){\small $\Sigma$}
\put(-167,40){\small $\Delta$}
\put(-179,58){\small $p$}
\put(-179,35){\small $q$}
\put(-102,10){\small $\Sigma$}
\caption{Producing an immersed disk $\Sigma$ in $X-\q$ bounded by $\sigma$}
\label{finger}}


If the {\sl geometric} intersection number of $\Sigma$ with the core $E$ of some 2-handle (say in $A_k$) is still nonzero, then choose a cancelling pair $p,q$ of points in $\Sigma\cap E$ and associated Whitney circle $\delta$.  This circle bounds an immersed disk $\Delta$ in the complement of $A_k$, since $\pi_1(X-A_k)=1$ (see \figref{finger}b again).  Pushing the core $F$ of each 2-handle that intersects $\Delta$ across $E$ by a finger move, we can arrange for $\Delta$ to lie in the complement of $\q$ (see \figref{finger}c). Here we are moving the corks rather than $\Sigma$, introducing extra intersection points between the cores of their $2$-handles.  This modifies $L$ by adding further clasps through new dotted circles in $J$, so does not alter the property that each $A_i$ is encased in $\q$. 
Now, a trivialization of the normal bundle of $\Delta$ induces a framing of $\delta$ which can be made to match the Whitney framing by ``boundary twisting"  $\Delta$ around $\Sigma$, at the cost of introducing additional intersection points between $\Delta$ and $\Sigma$ (see Freedman-Quinn \cite[\S1.3--1.4]{freedman-quinn:4-manifolds}).  The result is an immersed Whitney disk $\Delta$ for $\delta$, giving rise to a regular homotopy of $\Sigma$ that eliminates the intersection points $p$ and $q$, without adding any additional intersections between $\Sigma$ and $\q$ (see \figref{finger}d).  Repeating this procedure removes all the intersections of $\Sigma$ with $\q$.  

Since finger moves are supported near arcs, which by general position can be made to miss any finite collection of disks in $X$, this process can be iterated to produce null-homotopies for {\sl all} the belt circles of the 2-handles in $\q$, completing the proof.
\end{proof}


\vskip-5pt
\vskip-5pt

\head{Proof of the Consolidation \thmref{consolidation}}

\def\Kh{\widehat{K}}
\def\Lh{\widehat{L}}

Isotop the $A_i$ (for $i=1,\dots,n-1$) by the Finger \lemref{fingermove} so that each is encased in $\q$, which (as in the proof of the lemma) is equipped with an \ac\,structure extending ones $[K_i,L_i]$ on the $A_i$.  
The Encasement \lemref{encasement} extends this structure to one $[K,L]$ on an \ac\,cork $C_0$ in $X$ containing $\q$, and thus also encasing each $A_i$.   Let $C_i= (C_0)_{A_i, \tau_i}$, which is \ac\ by \remref{twistAC}.  By definition, the cork replacements $X_{C_0, C_i}$ are diffeomorphic to the cork twists $X_{A_i,\tau_i}$ rel boundary.

We argue that the multicork 
$
\calc \ = \ (C_0, C_1, \dots, C_{n-1})
$
is simple, that is, $C_0\cup -C_i \cong S^4$ for all $i$. 

First consider the double $C_0\cup -C_0$.  It has a handle structure extending $[K,L]$, obtained by attaching ``dual" $2$-handles along the belt circles $L^*$ of the $2$-handles of $C_0$ (i.e.\ the meridians of $L$) and ``dual" $3$-handles along the belt $2$-spheres of the $1$-handles of $C_0$, and then capping off with a $4$-ball. Hence the complement $(C_0\cup -C_0)-\interior(A_i)$ can be built from $\partial A_i \times [0,\varepsilon]$ by attaching the $1$-handles $K-K_i$ and $2$-handles $(L-L_i) \cup L^*$, plus some $3$-handles and a single $4$-handle.  Now since the $2$-handles $L- L_i$ homotopically cancel the $1$-handles $K-K_i$, each $2$-handle $J$ in $L-L_i$ can be slid over the dual $2$-handles $L^*$ 
 until it geometrically cancels its corresponding $1$-handle in $K-K_i$. The dual $2$-handle $J^*$ is then free to be cancelled with a $3$-handle. Proceeding inductively, we can cancel all the $2$-handles in $L^*-L^*_i$ with $3$-handles.  Since $A_i=[K_i,L_i]$ has not moved during this process, the handlebody structure for $(C_0\cup -C_0)-\interior(A_i)$ now consists solely of the dual $2$-handles attached along $L^*_i$ and the $3$-handles dual to the $1$-handles $K_i$. But, this is exactly $-A_i$. Thus each 
$
 C_i \cup -C_0 \ = \ (C_0\cup -C_0)-\interior(A_i) \,\cup_{\tau_i} -A_i
$ 
is diffeomorphic to the twisted double $A_i \cup_{\tau_i} -A_i$, which is in turn diffeomorphic to the $4$-sphere since $A_i$ is simple, and so $C_0\cup-C_i \cong S^4$ for all $i$, as asserted. 

Since $\calc$ is simple, its pinwheel $P$ (which is \ac\ by \remref{pwprop}) has an embedding in $X$ correlated with the embedding $C_0 \subset X$ (by \lemref{pinwheel}). Note that if $X$ is closed and each $A_i$ has simply-connected complement, the Finger \lemref{fingermove} gives an isotopy of the $A_i$ so that $\q$ has simply-connected complement. Hence, the embedding $C_0 \subset X$ and therefore the embedding $P \subset X$ can be made tight (by the Encasement and Pinwheel Lemmas \ref{encasement} and \ref{pinwheel}). Since $X$ is closed, this implies that $X-P$ is simply-connected. Thus $P$ is an \ac\ cork of order $n$ with the desired properties.
\qed 
\smallskip

\section{Main Results}\label{infinite}

The Consolidation Theorem implies the following version of the Finite Cork Theorem, extending the one stated in the introduction that treats the closed case.

\begin{namedtheorem}{Finite Cork Theorem} \label{rfct}
Let $X_i$ $(i\in\bz_n)$ be any finite list of compact simply-connected  $4\text{-manifolds}$ all homeomorphic to a given one $X = X_0$ that is closed or bounded by a homology sphere.   Then there is an \ac\ cork $(C,h)$ of order $n$ in $X$, with simply-connected complement in the closed case, whose twists $X_{C,h^i}$ are diffeomorphic to $X_i$ for each $i\in\bz_n$.  In the bounded case, these diffeomorphisms can be chosen to extend any given boundary identifications $f_i\colon\del X \to \del X_i$. 
\end{namedtheorem} 

\vskip-5pt
\vskip-5pt

\proof The closed case follows from the bounded case by removing the interior of a $4$-ball and noting that any orientation preserving diffeomorphism of the $3$-sphere is isotopic to the identity.  In the bounded case, the Relative Involutory Cork Theorem \ref{ICT} provides simple involutory corks $(A_1, \tau_1), \dots, (A_{n-1}, \tau_{n-1})$ with tight embeddings $A_i \subset X$ such that $f_i$ extends to a diffeomorphism $X_{A_i,h_i} \lto X_i$ for each $i$. The result is now immediate from the Consolidation Theorem \ref{consolidation}. \qed

\begin{remarks}
\ia\ When $X$ is closed, the cork in the Finite Cork Theorem and its complement can both be made Stein (or PC, borrowing the terminology from 
\cite{akbulut-matveyev:positron}) by extending the proof of Theorem 5 in \cite{akbulut-matveyev:positron} to the finite order case. In outline, first the complement of the multicork from the proof of the Consolidation Theorem \ref{consolidation} is made Stein, followed by each component of the multicork. This involves adding only homotopically cancelling $1,2$-handle pairs to each manifold. So, the multicork remains simple, and its complement simply-connected. Stein structures for the cork and its complement are then obtained by arguing that the pinwheel of an embedded simple order $n$ multicork with Stein components and complement is Stein (as shown in Part (3) of the proof of Theorem 5 in \cite{akbulut-matveyev:positron} when $n=2$) as is the complement of the correlated embedding.

\ib\ The Finite Cork Theorem allows one to generalize many conclusions that can be drawn from the Involutory Cork Theorem. The relative version in particular helps to \emph{localize} corks. For instance, suppose $E_i$ (for $0<i<n$) is any list of $4$-manifolds homeomorphic to some minimal elliptic surface $E = E(m)$, and obtained from it by log transforms along regular fibers.  These fibers 
can be  taken inside a Gompf nucleus $N$ of $E$ by \cite[Prop. 3.4]{gompf:nuclei}, 
so we obtain a list of manifolds $N_i\subset E_i$ homeomorphic to $N$, with diffeomorphic complements $E_i-N_i$. Since $N$ is simply-connected with homology sphere boundary (see \cite{gompf:nuclei}, Prop. 3.1),
the Finite Cork Theorem produces a cork $(C , h)$ of order $n$ in $N$ with $N_{C,h^i}$ is diffeomorphic to $N_i$ rel boundary for each $i$. Since all self-diffeomorphisms of $\partial N$ extend over $E - N$, by \cite[Lemma 3.7]{gompf:nuclei}, the cork twist $E_{C,h^i}$ is diffeomorphic to $E_i$ for each $i$.  Thus the $4$-manifolds $E_i$ are all obtained by twisting a single order $n$ cork inside a nucleus in $E$. 
\end{remarks}

The Relative Involutory Cork Theorem \ref{ICT} and the Consolidation Theorem \ref{consolidation} can also be used to pull apart finite families of corks in closed simply-connected  $4$-manifolds, in the following sense:

\begin{namedtheorem}{Separation Theorem}\label{separation} 
For any $n \in \mathbb{N}$ and any family of corks $(C_1,\sigma_1), \dots, (C_n,\sigma_n)$ embedded in a {\sl closed} simply-connected $4$-manifold $X$, there is a corresponding family of simple involutory corks $(A_i, \tau_i)$, embedded disjointly and each with simply-connected complement in $X$, whose twists $X_{A_i, \tau_i}$ are diffeomorphic to $ X_{C_i,\sigma_i}$ for each $i$. 
\end{namedtheorem}  

\proof 
The Involutory Cork Theorem gives the cork $(A_1,\tau_1)$.  Suppose inductively that corks $(A_1,\tau_1),\dots,(A_{n-1},\tau_{n-1})$ in $X$ have been found for some $n>1$ satisfying the theorem, so in particular $A_1,\dots,A_{n-1}$ are disjoint.  
By the Consolidation \thmref{consolidation}, there is a cork $(C, h)$ of order $n$ embedded with simply-connected complement in $X$ containing 
$$
A_1\cup\cdots\cup A_{n-1}\cup C_n
$$
whose cork twists $X_{C, h^i}$ are diffeomorphic to $X_i$ for each $i \le n$.  Note that the proof of \ref{consolidation} in general requires an initial isotopy of all the corks, but in this case it is only necessary to move $C_n$ since the corks $A_1, \dots, A_{n-1}$ are already disjoint. 

The Relative Involutory Cork \thmref{ICT} now yields a simple involutory cork $(P, \tau)$ tightly embedded in $Y = X-\interior(C)$ with a diffeomorphism $Y_{P,\tau} \to Y$ that restricts to $h^n \colon \del Y_{P,\tau} \to \del Y$, where both boundaries are identified with $\del C = -\del Y$.  Extending by the identity across $C$ gives a diffeomorphism $X_{P,\tau} \to X_{C, h^n}$, which can be composed with a diffeomorphism $X_{C, h^n} \cong X_n$ as above to give one $X_{P,\tau}\cong X_n$. Furthermore, since $P\subset Y$ is tight, $\pi_1(\partial C)$ maps onto $\pi_1(Y-P)$ and so $X-P$ is simply-connected.  Now set $(A_n,\tau_n) = (P,\tau)$, and proceed by induction. 
\qed

\smallskip

\begin{remark*}
Our proof of the Separation Theorem fails when $X$ has nonempty boundary, as we cannot guarantee the existence of a cork $C$ such that $X-C$ is simply-connected, and so cannot apply the Relative Involutory Cork Theorem in the inductive step. 
\end{remark*}


\begin{thebibliography}{10}

\bibitem{akbulut:contractible}
S.~Akbulut, \emph{A fake compact contractible {$4$}-manifold}, J. Differential
  Geom. \textbf{33} (1991), no.~2, 335--356. 

\bibitem{akbulut-matveyev:positron}
S.~Akbulut and R.~Matveyev, \emph{A convex decomposition theorem for $4$-manifolds}, Internat. Math. Res. Notices \textbf{7} (1998), 371--381. 

\bibitem{akbulut-ruberman:absolute}
S.~Akbulut and D.~Ruberman, {\em Absolutely exotic contractible $4$-manifolds}.  Comm. Math. Helv., \textbf{91} (2016), no.~1, 1--19. 
  
\bibitem{akbulut-yasui:corks-plugs}
S.~Akbulut and K.~Yasui, \emph{Corks, plugs and exotic structures}, J.
  G\"okova Geom. Topol. GGT \textbf{2} (2008), 40--82. 

\bibitem{ac}
J.~J.~Andrews and M.~L.~Curtis, \emph{Free groups and handlebodies}, Proc. Amer. Math. Soc. \textbf{16} (1965), no.~2, 192--195. 

\bibitem{akmr:equivariant}
D.~ Auckly, H.~ J.~ Kim, P.~ Melvin, and D.~Ruberman,
{\em Equivariant corks}, Alg. and Geom.\ Top.
\textbf{17} (2017), 1771--1783.

\bibitem{boyer}
S.~Boyer, 
{\em Simply-connected {$4$}-manifolds with a given boundary},
Trans. Amer. Math. Soc. \textbf{298}, no.~1, 331--357.

\bibitem{casson}
A.~Casson, \emph{Three lectures on new infinite constructions in $4$-dimensional manifolds}, (notes prepared by L. Guillou), A la Recherche de la Topologie Perdue, Progress in Mathematics vol. 62, Birkh\"auser, 1986, pp. 201-244. 

\bibitem{cerf}
J.~Cerf, \emph{Sur les diff\'eomorphismes de la sph\'ere de dimension trois $(\Gamma_4 = 0)$}, Lecture Notes in Math. 53, Springer-Verlag, Berlin (1968).

\bibitem{curtis-freedman-hsiang-stong}
C.~L.~Curtis, M.~H.~Freedman, W.~C.~Hsiang, and R.~Stong,
  \emph{A decomposition theorem for {$h$}-cobordant smooth simply-connected
  compact {$4$}-manifolds}, Invent. Math. \textbf{123} (1996), no.~2, 343--348.


\bibitem{freedman:simply-connected}
M.~H.~Freedman, \emph{The topology of four--dimensional manifolds}, J.
  Diff.\ Geo. \textbf{17} (1982), 357--432.

\bibitem{freedman-quinn:4-manifolds}
M.~H.~Freedman and F.~Quinn, \emph{Topology of {$4$}-manifolds}, Princeton University Press, Princeton, 1990.

\bibitem{gluck:2-spheres}
H.~Gluck, \emph{The embedding of two-spheres in the four-sphere}, Bull. Amer. Math. Soc. {\bf 67} (1961), 586--589. 

\bibitem{gompf:nuclei}
R.~ Gompf, \emph{Nuclei of elliptic surfaces}, Topology {\bf 30} (1991), no.~3, 479--511. 

\bibitem{kirby:4-manifolds}
R.~C.~Kirby, \emph{The topology of {$4$}-manifolds}, Lecture Notes in Mathematics, vol. 1374, Springer-Verlag, Berlin, 1989. 

\bibitem{kirby:cork}
\bysame, \emph{Akbulut's corks and h-cobordisms of smooth, simply-connected $4$-manifolds}, Turkish J. Math. \textbf{20} (1996), 85--93.

\bibitem{matveyev}
R.~Matveyev, \emph{A decomposition of smooth simply-connected
  {$h$}-cobordant {$4$}-manifolds}, J. Differential Geom. \textbf{44} (1996),
  no.~3, 571--582. 

\bibitem{stong}
R.~Stong, \emph{A structure theorem and a splitting theorem for simply-connected $4$-manifolds}, Math.\ Res.\ Letters \textbf{2} (1995), 497--503.

\bibitem{tange:finitecorks}
M.~Tange,
{\em Finite order corks}, 
Internat. J. Math. \textbf{28} (2017), 26.

\bibitem{tange}
\bysame,
{\em Non-existence theorems on infinite order corks}, preprint
\newblock \url{http://arXiv:1609.04344}, 2016.

\bibitem{trace}
B.~Trace,
{\em A class of $4$-manifolds which have $2$-spines}, Proc. Amer. Math. Soc. \textbf{79} (1980), 155--156.

\bibitem{wall:diffeomorphisms}
C.~T.~C. Wall, \emph{Diffeomorphisms of {$4$}-manifolds}, J. London Math. Soc. \textbf{39} (1964), 131--140. 

\bibitem{wall:4-manifolds}
\bysame, \emph{On simply-connected {$4$}-manifolds}, J. London Math. Soc. \textbf{39} (1964), 141--149. 

\bibitem{yasui}
K.~Yasui,
{\em Nonexistence of twists generating exotic $4$-manifolds}, Trans. Amer. Math. Soc. \textbf{372} (2019), 5375--5392.


\end{thebibliography}
\end{document}